\documentstyle[10pt,amscd]{amsart}

\subjclass{55N10, 55N45, 55R20, 55U99}
\keywords{Simplicial set, cubical set,  fibration,
    $K(\pi,n)$-space, cross section, twisted tensor product, 
    perturbed differential, algebraic 
    model, spectral sequence}

\newtheorem{lemma}{\ \ \ Lemma}[section] 
\newtheorem{theorem}{\ \ \ Theorem}[section]
\newtheorem{corollary}{\ \ \ Corollary}[section]
\newtheorem{proposition}{\ \ \ Proposition}[section]

\theoremstyle{definition}
\newtheorem{definition}{\ \ \ Definition}[section]

\numberwithin{equation}{section}

\begin{document}

\author{N. Berikashvili}

\title[An algebraic model of fibration]
    {An algebraic model of fibration with 
    the fiber $K(\pi,n)$-space}

\begin{abstract}
For a fibration with the fiber $K(\pi,n)$-space, the algebraic model as
a twisted  tensor product of chains  of the base with standard chains
of $K(\pi,n)$-complex is given which preserves multiplicative structure
as well. In terms of this model   the action of the $n$-cohomology
of the base with coefficients in  $\pi$ on the homology of fibration is 
described. 
\end{abstract}
\maketitle 
\section{Introduction}

For a fibration $F @>>> E @>>> B$ there is the Brown model [1], [2] as a twisted
tensor product $C_{*}(B) \otimes_{\varphi} C_{*}(F)$. However, this model
gives us no information about the multiplicative structure. The aim of this 
paper is to construct a model for the particular case of the fiber 
$F=K(\pi,n)$ which would inherit  the multiplicative structure of $C^{*}(E)$ 
as well. The model is given as a twisted tensor product 
 $C_{*}(B) \otimes_{z^{n+1}} C_{*}^{\Box}(L(\pi,n))$
of the singular chain complex of the base with the chain complex of the cubical version
of the complex $K(\pi,n)$ [3], [4] and it describes the cubical singular chain 
complex of the total space (Theorems 5.1 and 6.1 below). It turns out
that this model is actually the  chain complex of a cubical complex and hence, 
in addition, it carries a
Serre cup product structure as well. In Section 7 the action of the
group $H^{n}(B,\pi)$ on the model is discussed.

This paper is a reply to the needs of obstruction theory and is applied
in [5]. The main result was announced in [6].

\section{Preliminaries}

Let $X$ be a $CW$-complex and $F^{p}X= X^{p}$ be its filtration by  skeletons.
If $H$ is the singular homology theory with coefficient
group $G$, then for 
the first term of the related spectral sequence we have $E_{p,q}^{1}=0,\;\; q > 0,$ 
and $(E_{p,0}^{1}, d^{1})$ is the chain complex
$$  (H_{p}(X^{p},X^{p-1},G),d^{1})     $$  
with the homology isomorphic to the singular homology of the space $X$. 
Of course $H_{p}(X^{p},X^{p-1},G)$ is isomorphic to $\Sigma G_{\sigma}$, 
where $\sigma$ are $p$-cells of $X$ and $G_{\sigma}=G$. The above chain complex
will be referred to as the cell chain complex of the $CW$-complex.

In some cases the chain map of the cell-chain complex of the $CW$-complex 
to the singular or the cubical singular chain complex
is defined 
which  induces the sisomorphism of homology  stated above.

Recall that the cubical set  is a sequence of 
sets $Q=\{ Q_{0},Q_{1},Q_{2} \cdots \}$ together with the boundary and 
degeneracy operators
\begin{align*}
    d_{i}^{\epsilon} & :  Q_{n} \rightarrow Q_{n-1}, \;\;\; 1 \leq i \leq n, 
        \;\;\; \epsilon = 0,1,  \\
    s_{i} &:  Q_{n} \rightarrow   Q_{n+1},\;\;\; 1 \leq  i \leq  n+1,
\end{align*}
subject to the standard equalities (see, e.g., [4], [7]).

$Q_{n}$ are  $ n$-cubes of $Q; \; \; \; 
\sigma^{n} \in Q_{n}$ is said to be degenerate if 
$$  \sigma^n = s_{i} \tau^{n-1}.        $$
The abbreviation for $d_{i}^{\epsilon} \sigma^n$ is 
$\sigma^{n,\epsilon}_i,\;\;\epsilon =0,1.$

The main example of a cubical set is  the cubical singular complex $ Q(X)$
of a space $X.$

Milnor's notion of the realization of a simplicial set [8] works for  
a cubical set as well and runs as follows.

Let $I^{n}$ be the standard $n$-cube and let
\begin{gather*}
    e_{i}^{\epsilon}: I^{n-1} @>>> I^{n} ,\;\;\; 1 \leq i \leq n , \;\; 
         \epsilon =0,1, \\
    p_{i}: I^{n+1} \rightarrow I^{n}, \;\;\; 1 \leq i \leq (n+1),    
\end{gather*}
be $i$-face imbeddings and $i$-projections.

Let Q  be a cubical set. Then
the realization $|Q|$ is defined as a factor set of the space
$$  \bigcup_n  (Q_{n} \times I^{n})     $$
by the identifications
\begin{gather*}
    (\sigma^{n}, e_{i}^{\epsilon}x)=(d_{i}^{\epsilon}(\sigma^{n}),x), 
        \;\;\;\sigma \in Q_{n}, \;\; x \in I^{n-1},  \\
    (s_{i}(\sigma^{n}), x)=(\sigma^{n},p_{i}x) ,
        \;\;\; \sigma \in Q_{n}, \;\;\; x \in I^{n+1}, \;\;    
        1 \leq i \leq (n+1).   
\end{gather*}
$|Q|$ is a CW-complex with $n$-cells in  1--1 correspondence with the
nondegenerate $n$-cubes of $Q$.

There is a standard continuous map
$$  r:|Q(X)| @>>> X.    $$

We define the Giever-Hu realization [9], [10] of the cubical set $Q$ and denote 
it by $||Q||$, omitting in the above definition the equality with degeneracy 
operators. Then $||Q||$ has cells in 1--1 correspondence with all 
cubes of $Q.$ We have the canonical continuous map
$$  ||Q|| @>>> |Q|.     $$

We recall the notion  of homology for the cubical set $Q$. 
$\tilde{C}_{*}^{\Box}(Q)$ is the chain complex spanned in each dimension $n$
by all $n$-cubes of $Q$, and the boundary operator $d$ is defined by
$$  d\sigma^{p}= \Sigma(-1)^{i}
    d^{0}_{i}\sigma^{p}-\Sigma(-1)^{i} d^{1}_{i} \sigma^{p}.    $$
$\overline{C}\,{}_{*}^{\Box}(Q)$
 is a chain subcomplex of
$\tilde{C}_{*}^{\Box}(Q)$ 
spanned by the degenerate cubes. 
The quotient complex
$$  C_{*}^{\Box}(Q)=\tilde{C}_{*}^{\Box}(Q)/
        \overline{C}\,{}_{*}^{\Box}(Q)      $$ 
is said to be the chain complex of $Q$ and its 
homology is called the homology of $Q$.

If $G$ is a group of coefficients, then by 
$$  C^{n}_{\Box}(Q,G)   $$
we understand the normalized cochains of $Q$, i.e., those which are 
zero on the degenerate cubes.

For the singular cubical set $Q(|Q|)$ of the space $|Q|$ consider the standard 
imbedding
$$  Q \subset Q(|Q|).   $$
It induces the following isomorphisms of homology and cohomology:
$$  H_{n}^{\Box}(Q,G)=H_{n}^{\Box}(|Q|,G),    \;\;\;   
      H^{n}_{\Box}(Q,G)=H^{n}_{\Box}(|Q|,G).        $$
Hence we find that the map
$ r: |Q(X)| \rightarrow X$ induces the isomorphisms 
$$  H_{n}^{\Box}(|Q(X)|,G)=H_{n}^{\Box}(X,G)),\;\;\;
        H^{n}_{\Box}(|Q(X)|,G)=H^{n}_{\Box}(X,G)).          $$

The interaction of simplicial and cubical sets is as follows.
Let $\Delta_{n}$ be the standard simplex and 
consider the map
\begin{equation}
    \psi_{n} : I^{n} @>>> \Delta_{n} 
\end{equation}
from [11] defined by
\begin{equation}
    \begin{array}{c}
        y_{0}=1-x_{1},    \\
        y_{1}=x_{1}(1-x_{2}),   \\
        y_{2}=x_{1} x_{2}(1-x_{3}),     \\
        \vdots  \\
        y_{n-1}=x_{1} x_{2}x_{3} \cdots x_{n-1}(1-x_{n}),   \\
        y_{n}=x_{1} x_{2}x_{3}  \cdots x_{n}.  
    \end{array}
\end{equation}

Clearly, $\psi_{n}$ is the map of the pairs 
$ \psi_{n} :(I^{n}, \partial I^{n} ) @>>> ( \Delta_{n},\partial \Delta_{n}), $
inducing the homeomorphism
$I^{n}/ \partial I^{n} @>>> \Delta_{n}/ \partial \Delta_{n}.$

In [11], for a singular simplex 
$\sigma^{n} :\Delta_{n} @>>> X$  J. P. Serre considers the
composition $I^{n} @>\psi_{n}>> \Delta_{n} @>\sigma^{n}>> X $ as the singular 
cube $t\sigma^{n}$
of $X$, defining the chain  map 
$C_{*}(S(X)) @>>> C_{*}^{\Box}(Q(X))$ from the simplicial 
singular 
chain complex to the singular normalized cubical chain complex
(see the identities below).
Hence for an abelian group $G$ we have the cochain map
\begin{equation}
    C^{*}_{\Box}(Q(X),G) @>>>C^{*} (X,G)    
\end{equation}
from  the normalized singular cubical 
cochain complex to an ordinary singular cochain complex.

For map (2.3) at the geometrical level, $ t:S(X) @>>> Q(X) $, we have
 the identities
\begin{equation}
\begin{gathered}
    d^{1}_{i}t(\sigma^{n})=
        t(\sigma^{n}_{i-1}),\;\;\; 1 \leq i \leq n,  \\
    d^{0}_{n}t(\sigma^{n})=t(\sigma^{n}_{n}),  \\
    d^{0}_{i}t(\sigma^{n})=
        \underbrace{s_{i}s_{i} \cdots s_{i}}_{n-i} 
        t(\sigma^{n}_{n,(n-1),\cdots ,i}), \;\;\; 1 \leq i <n.
\end{gathered}
\end{equation}
 We see that $(t\sigma)_{i}^{0}$ is a degenerate cube 
 if $i \neq n.$

The above equalities imply that the following lemma is valid.

 \begin{lemma}
Let $0 < k \leq n$ and let $P$ be the sum of closed $k$-cubes of $I^{n}$ 
which are not degenerated by $\psi_{n}$. Then $P$ is a deformation 
retract of $I^{n}.$
\end{lemma}

One defines a chain map
$C_{*}^{\Box}(X) @>>> C_{*}^{N}(X)$
of the normalized complexes as follows. Consider the standard
triangulation  of $I^{n}$ and let $s(I^{n})$ 
denote the basic $n$-dimensional chain of this triangulation as a singular 
chain of the space $I^{n}$. Define $f$ by
\begin{equation}
    f(\sigma^{n})= \sigma^{n}_{*} (s(I^n)) \in C_{n}(X). 
\end{equation}
If $\sigma^n$ is a degenerate cube, then all $n$-simplexes in  $f(s(I^n))$ 
are degenerate and therefore $f $ defines the map of the normalized complexes.

For cubical cochains  J. P. Serre defines the $\smile$-product which  
is certainly valid for general cubical sets too [11].
The definition runs as follows.
For the subset   
$$  K=(  i_{1} < i_{2} <i_{3}  \cdots  <i_{k} ) 
        \subset (1,2,3, \cdots, n)      $$
and   $ \sigma^{n} \in Q_{n}$ we introduce the notation  
$$  d_{K}^{\epsilon} \sigma^{n} =d_{i_{1}}^{\epsilon}  
         d_{i_{2}}^{\epsilon} d_{i_{3}}^{\epsilon}
        \cdots d_{i_{k}}^{\epsilon}(\sigma^{n}),\;\; \; 
        \epsilon =0,1.      $$
 
 Let $c^{p} \in C^{p}(Q,\Lambda),\;\;\;    c^{q} \in C^{q}(Q,\Lambda)$ 
 be the normalized  $p$- and 
$q$-cochains of $Q$ with coefficients in  a commutative ring $\Lambda$.
 Define the product $c^{p+q}=c^{p} \smile c^{q}$     by
\begin{equation}
    c^{p+q}(\tau^{p+q})=\Sigma_{(H,K)}(-1)^{a(H,K)}c^{p}
    (  d_{K}^{0}(\tau^{p+q})) c^{q}(d_{H}^{1}(\tau^{p+q})), 
\end{equation}
where $(H,K)$ is the decomposition  of  
 $\{1,2,3, \cdots,  (p+q) \} $  into two disjoint subsets.

Note that if $\tau^{p+q}$ is of the form $t(\sigma)$, then on 
the right  side only one summand, precisely that of  the decomposition 
$(1,2,3, \cdots,p) \bigcup  (p+1,p+2, \cdots ,p+q),$ is not zero, which 
shows that the cochain map (2.3) preserves multiplicative 
structure.  

For any abelian group $\pi$  and any positive  integer $n$ one introduces 
the  cubical version of the $ K(\pi, n)$-complex [3], [4].  
$ L(\pi, n)$ is a cubical set with $p$-cubes being $n$-cocycles, 
$Z^{n}(I^{p},\pi)$,
where $I^{p} =I \times I \times \cdots \times I$ is the standard cube with  
standard faces, and the cochain is understood as a cell cochain. The boundary 
and degeneracy operators  are defined similarly to those in  the c.s.s. case.

For the complex $L (\pi, n)$  we have, in  $C_{*}^{\Box}(L (\pi, n))$,  
the  ring structure that converts  $ C_{*}^{\Box}(L (\pi, n))$ into a graded 
differential  ring and hence $H_{*}^{\Box}(L(\pi,n))$ into a graded 
 ring. The multiplication is defined as 
follows. Given cubes $\sigma^{p}$ and $\sigma^{q} $ in $L (\pi, n)$, i.e.,  
$\sigma^{p} \in Z^{n}(I^{p}, \pi)$, $\sigma^{q} \in Z^{n}(I^{q}, \pi)$, we 
define the $(p+q)$-cube $\sigma^{p} \circ \sigma^{q}= 
\sigma^{p+q} \in Z^{n}(I^{p+q}, \pi) $ as 
$pr^*_{1} \sigma^{p} + pr^*_{2} \sigma^{q}$, where 
$$  pr_{1}:I^{p+q} \rightarrow I^{p}, \; \; \;    
        pr_{2}:I^{p+q}\rightarrow I^{q}     $$
are standard projections of $I^{p+q}$ on its first $p$-cube and its last 
$q$-cube.
Under this multiplication $\tilde{C}_{*}^{\Box}(L(\pi,n))$ is the differential 
graded ring and $\bar{C}_{*}^{\Box}(L(\pi,n))$ is the ideal, and therefore  
$\tilde{C}_{*}^{\Box}(L(\pi,n))/\overline{C}\,{}_{*}^{\Box}(L(\pi,n))  $ 
is the graded differential associative ring. We have  
$
C_{0}^{\Box}(L(\pi,n))\!=\!Z$, $C_{i}^{\Box}(L(\pi,n))
   =0$, $0<i<n$, while  
$1 \in  Z$ is the unitary element of the algebra $C_{*}^{\Box}(L(\pi,n)).$

The Eilenberg--MacLane complex $K(\pi,n)$ is a simplicial set and  \linebreak
$C_* (K(\pi,n))$ is the skew-commutative differential ring. The above map
$$  \psi_n : I^{n} @>>> \Delta_{n}      $$ 
enables us to define the chain map
$$  C_* (K(\pi,n)) @>>>C_{*}^{\Box}(L(\pi,n))   $$   
similarly to that of J. P. Serre. In an obvious 
manner similar to that used for (2.5) we  define  
$$  C_{*}^{\Box}(L(\pi,n)) @>>>C_{*}^{N}(K(\pi,n))  $$
which is a homomorphism of differential rings. Hence $H_{*}^{\Box}(L(\pi,n))$ 
and $H_{*}(K(\pi,n))$ are isomorphic as rings and 
thus $H_{*}^{\Box}(L(\pi,n))$ is a commutative ring.

Recall some notions from the theory of perturbed differentials [12]. Let $Y$ be 
a filtered left graded differential module over a filtered graded differential 
algebra $ A$ under the pairing 
$$  A \otimes Y @>>> Y.         $$
Let $F^{i}A \supset F^{i+1}A$ and $F^{i}Y  \supset F^{i+1}Y$ be the 
above-mentioned filtrations.

Fix the integer $n$.
Let 
$$  T^n (A)=\{ a \in F^n A, \; |a|=+1, \;  d_{A}a=a a\};    $$
$a$ is known as a twisting element of the differential algebra. 

For $a \in T^n (A) $
the perturbed 
differential in $Y$ is defined as
$$  d_{a}x=d_{Y}x+ax, \; \; \; x \in Y.         $$
The filtered complex with this perturbed differential  will be denoted by 
$Y_{a}$. It is a 1-graded and filtered differential chain complex.
 Consider the set   $G^n (A)$   of   all elements  in  $A$ of the form
$$  g=1+p, \; p \in F^1 A, \; d_A p \in F^n A.      $$
It is obvious that this set is the group under the product operation 
in the algebra $A$.
The group $G^n(A)$ acts from the left on the set 
$T^n (A)$ by 
\begin{equation}
    \bar{a} =(1+p)\circ a=(1+p)^{-1}a (1+p)+(1+p)^{-1}d_{A}(1+p)
\end{equation} 
which is equivalent to the equality
\begin{equation}
    d_{A}(1+p)=(1+p)\bar{a}-a(1+p).  
\end{equation}
The set of orbits under action (2.7) will be  denoted by  $D^n (A)$.

For 
$\overline{a}=(1+p)\circ a$
consider the map
$\varphi _{p}: Y_{\bar{a}}\rightarrow Y_{a} $ 
given by 
\begin{equation}
    \varphi_{p}(x) =(1+p)x=x+px.    
\end{equation}
This is a chain map and, obviously,  1--1.

\begin{proposition}
 If $f:A_1 @>>> A$ is a morphism of filtered differential
algebras inducing an isomorphism of the $n$th terms of spectral sequences,   
then the induced map $D^nf : D^n( A_1 ) @>>> D^n(A)$ is $1$--$1$ $[12]$.
\end{proposition}

\section{Auxiliary  Complexes}

Let $B$ be a topological space  and $\pi$  an abelian group. 
Let  $L(\pi ,n)$ be the $K(\pi,n)$-complex of the preceding section and 
$C_{*}^{\Box}(L(\pi ,n))$ its normalized chain  complex.  Let $G$ be an abelian 
group. Consider the chain and cochain complexes $C_{*}^{\Box}(L(\pi ,n),G)$ 
and $C^{*}_{\Box}(L(\pi ,n),G)$.

Consider the complexes
\begin{equation}
\begin{aligned}
    Y_{**}(B,\pi,n) & =C_{*}(B,C_{*}^{\Box}(L(\pi ,n))),  \\
    Y_{**}(B,\pi,n,G) & =C_{*}(B,C_{*}^{\Box}(L(\pi ,n),G)),   
\end{aligned}
\end{equation}
i.e., the group of singular chains of $B$ with coefficients in the normalized 
integral chains of the complex $L(\pi,n)$ and the group of singular chains 
of $B$ with coefficients in the  normalized chains of the complex 
$L(\pi,n)$ with coefficients in $G$. These complexes are bigraded differential 
modules, the first differential being that  of $B$ and the second one that 
of $L(\pi,n)$. Both complexes are covariant functors on the category of 
topological spaces. It is clear that   
\begin{align*}
    Y_{**}(B,\pi,n) & =C_{*}(B) \otimes C_{*}^{\Box}(L(\pi ,n)),  \\
    Y_{**}(B,\pi,n,G) & =C_{*}(B) \otimes C_{*}^{\Box}(L(\pi ,n),G),  \\
    Y_{**}(B,\pi,n,G) & =C_{*}(B,C_{*}^{\Box}(L(\pi ,n))) \otimes G. 
\end{align*}
The cohomology version is
$$  Y^{**}(B,\pi,n,G)=C^{*}(B,C^{*}_{\Box}(L(\pi ,n),G)).       $$
We have  
\begin{align*}
    Y^{**}(B,\pi,n,G) & =Hom(C_{*}(B), C^{*}_{\Box}(L(\pi ,n),G)),  \\
    Y^{**}(B,\pi,n,G) & =Hom(Y_{**}(B,\pi,n),G)).
\end{align*}

 Consider the graded differential algebra
$C_{*}^{\Box}(L(\pi,n))$ of the preceding section and the bigraded 
differential algebra
\begin{equation}
    A^{*}{}_{*}(B,\pi,n)=C^{*}(B,C_{*}^{\Box}(L(\pi,n))) 
\end{equation}
of  singular cochains of $B$ with coeficients in the normalized chains
of $L(\pi,n)$, where the multiplication is defined as the $\smile$-product of 
$B$ when the coefficients are multiplied by the product in 
the algebra $C_{*}^{\Box}(L(\pi,n))$. We obtain 
$$  A^{*}{}_{ *}(B,\pi,n)=Hom(C_{*}(B),C_{*}^{\Box}(L(\pi,n)))      $$
and the above multiplication is the same as that  defined by the composition 
$$  C_{*}(B) @>\Delta>>  C_{*}(B) \otimes  C_{*}(B) @>x\otimes y>>
     C_{*}^{\Box}(L(\pi,n))\otimes C_{*}^{\Box}(L(\pi,n)) @>\mu>> C_{*}^{\Box}(L(\pi,n))    $$
for $x,y \in Hom(C_{*}(B),C_{*}^{\Box}(L(\pi,n))) $, 
where $\Delta$ is the coproduct  of $B$ and $\mu $ is the product in 
$C_{*}^{\Box}(L(\pi,n)).$ By lifting the second index in (3.2) we obtain 
the bigraded algebra of the fourth quadrant with both differentials 
increasing the degree of elements by +1:

$$  A^{*,-*}(B,\pi,n)=C^{*}(B,C^{-*}_{\Box}(L(\pi,n))).     $$

The total complex is the direct product of bihomogeneous
components 
$$  A^m =\Pi^{p,-q}_{p-q=m}A^{p,-q}     $$
and therefore every element $x^n$ is given uniquely as the sum of its components
$$  x^n =x^{n,o} +x^{n+1,-1}+x^{n+2,-2}+x^{n+3,-3}+ \cdots\,.   $$
The unitary element of the algebra  $1 \in A^{0,0} \subset
A^0 $  is the cochain of  $ B$,  $c^{0,0}$, equal to $1$ of $C_{*}(L(\pi,n))$ for every
$0$-simplex 
of $ B$.
$A^{*,-*}(B,\pi,n)$ is a contravariant functor both on the category of 
topological spaces and on the category of simplicial sets. The filtration 
of bialgebra $A$ with respect to the first degree is the decreasing one and 
is complete in the Eilenberg-Moore sense. The second term of the spectral 
sequence is
$$  E^{p,-q}_{2} = H^{p}(B,H_{q}(\pi,n)),       $$ 
the spectral sequence converging to  
$$  H^{p-q}(A^{*,-*}(B,\pi,n))  $$ 
in the Eilenberg--Moore sense. 

$Y_{**}(B,\pi,n,G)=Y^{-*,-*}(B,\pi,n,G)$ and $Y^{**}(B,\pi,n,G)$ 
are the left modules 
over the bigraded algebra $A^{*,-*}(B,\pi,n)$ to be understood as follows. 
The pairing 
\begin{equation}
    C_{p}^{\Box}(L(\pi,n)) \otimes C_{q}^{\Box}(L(\pi,n))  @>>> C_{p+q}^{\Box}(L(\pi,n)) 
\end{equation}
induces the pairings
\begin{align}
    C_{p}^{\Box}(L(\pi,n)) \otimes C_{q}^{\Box}(L(\pi,n)) \otimes G 
        & @>>> C_{p+q}^{\Box}(L(\pi,n))\otimes G,  \\ 
    C_{p}^{\Box}(L(\pi,n)) \otimes C^{p+q}_{\Box}(L(\pi,n),G)  
        & @>>> C^{q}_{\Box}(L(\pi,n),G). 
\end{align}

The bicomplex $Y_{**}(B,\pi,n,G)$ is the module over the algebra 
$A^{*}{}_{ *}(B,\pi,n)$ via the $\frown$-product of $B$ and the product of 
coefficients (3.3) and (3.4). The bicomplex $Y^{**}(B,\pi,n,G)$ is 
the left module over the algebra $A^{*}{}_{ *}(B,\pi,n)$
via  the $\smile$-product in $B$ and the product of coefficients (3.5). 

The spectral sequence arguments show that if $f:B_1 @>>> B $ induces 
an isomorphism of homology, then the induced
homomorphisms of complexes 
$$  A^{*,-*}(B,\pi,n) @>>> A^{*,-*}(B_1 , \pi,n)    $$   
and  
$$  Y_{**}(B_{1},\pi,n)@>>>Y_{**}(B,\pi,n)  $$    
induce an isomorphism of homology.
  
In particular, if $K$ is the ordered simplicial complex then the imbedding 
$K \subset S(|K|)$ induces the isomorphism 
\begin{equation}
    H(A^{*,-*}(|K|,\pi,n)) @>>> H(A^{*,-*}(K,\pi,n)). 
\end{equation}

Consider the algebra $A^{*,-*}(B,\pi,n)$ as a filtered differential algebra 
with the filtration defined by the first degree and introduce the notation 
\linebreak  $D(A^{*,-*}(B,\pi,n))=D^{n+1}(A^{*,-*}(B,\pi,n)).$ 
Then
$$  D(A^{*,-*}(B,\pi,n))=T^{n+1}(A^{*,-*}(B,\pi,n))
        /G^{n+1}(A^{*,-*}(B,\pi,n)),        $$
where $T^{n+1}(A^{*,-*}\!(B,\pi,n))$ is the set of all twisting elements of
$A^{*,-*}\!(B,\pi,n)$ of the form 
$$  a=a^{n+1,-n}+a^{n+2,-n-1}+a^{n+3,-n-2} \cdots,      $$ 
and $G^{n+1}(A^{*,-*}(B,\pi,n))$ is
the set of all elements  in  $A^{*,-*}(B,\pi,n)$ of the form
$$      1+p=1+p^{n,-n}+p^{n+1,-n-1}+p^{n+2,-n-2}+ \cdots    $$
(by virtue of $C_{i}^{\Box}(L(\pi,n))=0, \;\; 0 < i < n$).

\begin{proposition}
If  $f: B_1 @>>>  B$ induces  an isomorphism  of homology then 
$D(f) :D(A^{*,-*}(B,\pi,n))@>>> D(A^{*,-*}(B_1,\pi,n))$ is the $1$--$1$ map.
\end{proposition}

\begin{pf}
We   have $ E^{p,-q}_{2}=H^{p}(B,H_{q}(\pi,n)).$ By assumption,  $f$ 
induces an isomorphism of the second terms and  hence of the $(n+1)$th terms 
of spectral sequences. To complete the proof apply Proposition 2.1.  
\end{pf}

\section{ Transformation of the Functors
 $\beta: D(A^{*,-*}(B,\pi,n)) \to H^{n+1} (B,\pi)$
 and $\gamma: H^{n+1} (B,\pi) \to D(A^{*,-*}(B,\pi,n))$}

Consider the homomorphism  $\alpha : C_{n}(L(\pi,n)) @>>> G $ defined on 
the $n$-cubes of $L(\pi,n), \; \; \; \tau^{n}=[g],\; g \in  G,$
by $\alpha ([g])=g.$ It induces  the homomorphism 
of cochain complexes   $\beta : C^{*}(B,C_{n}(L(\pi,n))) @>>> C^{*}(B,\pi)$.
The latter homomorphism defines the map of sets 
$\beta : A^{*,-*}(B,\pi,n) @>>> C^{*}(B,\pi)$ by
$\beta (a)= \beta (a^{n+k, -n})$ if $a=a^{n+k,-n}+a^{n+k+1,-n-1}+ \cdots .$
Let $a \in  T(A(B,\pi,n)).$ Then  $|a|=1,\; \; da=aa $, 
and hence $\delta _{B}a^{n+1,-n} +d_{L}a^{n+2,-n-1}=0$.  
We see  that $\beta (d_{L}a^{n+2,-n-1})=0$, so that 
$\beta(\delta _{B}a^{n+1,-n})=0$ and $\delta _{B} \beta (a^{n+1,-n})=0.$
Hence,  if $a \in T(A^{*,-*}(B,\pi,n))$, then $\beta (a) \in  Z^{n+1}(B,\pi).$
 
If $\bar{a} \sim a,$ then $\beta (\bar{a})- \beta (a) \in B^{n+1}(B,\pi).$ 
Indeed, by (2.8) we have $d_{A}p= \bar{a} - a + p \bar{a} - ap,$ and 
therefore 
 $\delta_{B}p^{n,-n} + d_{L}p^{n+1,-n-1}= \bar{a}^{n+1,n}-a^{n+1,-n}. $
Hence $\delta_{B} \beta (p^{n,-n})= 
\beta(\bar{a}^{n+1,-n}) - \beta(a^{n+1,-n})$ 
which is equivalent to  $\delta_{B}\beta(p)= \beta (\bar{a})-\beta (a).$
The above shows that  the map  
$$  \beta :  D(A^{*,-*}(B,\pi,n)) @>>> H^{n+1}(B,\pi)       $$
defined by
$\beta(a) \in \beta (d), d \in D(A^{*,-*}(B,\pi,n)), a \in d,$
is correct.
 
The second transformation 
$$  \gamma : H^{n+1}(B,\pi) @>>> D(A^{*,-*}(B,\pi,n))   $$
is more difficult to define, since it requires the construction of
the mappings 
\begin{align*}
    Z^{n+1}(B,\pi)  & \rightarrow T^{n+1}(A^{*,-*}(B,\pi,n)),  \\
    Z^{n+1}(B,\pi) \times   C^{n}(B,\pi) & \rightarrow 
        G^{n+1}(A^{*,-*}(B,\pi,n))  ,
\end{align*}
with suitable properties.

 \begin{lemma}    
There is a rule which assigns  a twisting cochain 
$a(z^{n+1}) \in A^{*,-*}(B,\pi,n)$ to every cocycle 
$ z^{n+1} \in Z^{n+1}(B,\pi)$. If $f:B_{1} \rightarrow B$ is a map, 
then $a(f^{*}(z^{n+1}))=f^{*}(a(z^{n+1}))$.
\end{lemma}

\begin{pf}
For  $ z^{n+1} \in Z^{n+1}(B,\pi)$ we first construct the map
\begin{equation}
    \kappa_{z^{n+1}}: S(B) @>>> L(\pi,n),   
\end{equation}
$\kappa_{z^{n+1}}(\sigma^{m+1})$ being the $m$-cube $\tau ^{m} \in L(\pi,n) $  
with the following properties:
 
 (i)   $ \kappa_{z^{n+1}}(\sigma^{m+1})$ if $(m+1) \leq n$ is the unique 
 m-cube of $L(\pi,n)$ (which is degenerate for $m>0$);
 
 (ii) 
$ \kappa_{z^{n+1}}(\sigma^{n+1})=[g] \in \tilde{C}_{n}^{\Box}(L(\pi,n)),$
 where  $g=z^{n+1}(\sigma^{n+1});$

(iii)\ \vskip-0.8cm
\begin{equation}
\begin{aligned}
    d^{1}_{i}\kappa_{z^{n+1}}(\sigma^{m+1}) & 
        =\kappa_{z^{n+1}}(\sigma^{m+1}_{i}),\; \; \; \; 
        1 \leq i \leq m,  \\
    d^{0}_{1}\kappa_{z^{n+1}}(\sigma^{m+1}) & 
        =\kappa_{z^{n+1}}(\sigma^{m+1}_{0}),  \\
    d^{0}_{m}\kappa_{z^{n+1}}(\sigma^{m+1}) & 
        =\kappa_{z^{n+1}}(\sigma^{m+1}_{m+1}),  \\
    d^{0}_{i}\kappa_{z^{n+1}}(\sigma^{m+1}) & =
        \kappa_{z^{n+1}}(\sigma^{m+1}_{(1)})
        \kappa_{z^{n+1}}(\sigma^{m+1}_{(2)}),
         \; \; 1<i<m,
\end{aligned}
\end{equation}
where $\sigma^{m+1}_{(1)}$ and $\sigma^{m+1}_{(2)}$ are respectively
the first $i$-dimensional and the last $(m+1-i)$-dimensional face of 
$\sigma^{m+1}.$ 

With such $\kappa_{z^{n+1}}$ constructed, the twisting element 
$\tilde{a}_{z}$ in the auxiliary algebra 
$\tilde{A}^{*}{}_{ *}(B,\pi,n)=C^{*}(B,\tilde{C}_{*}(L(\pi,n)))$ is defined as 
$\tilde{a}=\tilde{a}^{2,-1}+\tilde{a}^{3,-2}+ \cdots$, where 
$\tilde{a}^{m+1,-m}(\sigma^{m+1})= \kappa_{z^{m+1}}(\sigma^{m+1}) \; 
\in \tilde{C}_{m}(L(\pi,n)). $

The above equalities imply that $d \tilde{a}=\tilde{a} \tilde{a}$. 
The image of $\tilde{a}$, $a(z^{n+1})$, with respect to the homomorphism
$\tilde{A}^{*}{}_{ *}(B,\pi,n) @>>> A^{*}{}_{ *}(B,\pi,n)$ induced by 
$\tilde{C}_{*}^{\Box}(L(\pi,n)) @>>> C_{*}^{\Box}(L(\pi,n))$ is 
a twisting cochain. By virtue of the fact that all cubes of positive dimension 
$< n$ are degenerate, we find that 
$$  a(z^{n+1})=a^{n+1,-n}+a^{n+2,-n-1}+ \cdots  $$  
and
$$  \beta(a)=\beta(a^{n+1,-n})=z^{n+1}.     $$

The rest of the proof consists in constructing $\kappa_{z^{n+1}}$.

The geometric background for the (algebraic) definition of 
$\kappa_{z^{n+1}}$ is as follows. Consider a cross section 
$s$ of the fibration over the $n$-skeleton of $|S(B)|$ such that
the obstruction cocycle $z^{n+1}(s)$ is equal  
to  $z^{n+1}, \;\; z^{n+1}=z^{n+1}(s).$
For  $\sigma^{m+1} \in S(B)$ consider the map 
 $I^{m+1} @>\psi>> \Delta_{m+1} @>\sigma^{m+1}>> |S(B)|.$
 Consider  a subcomplex of $n$-cubes of the cubical 
complex $I^{m+1}$ which is sent by $\psi$ in the $n$-simplex  of $\Delta_{m+1}.$ 
By virtue of Lemma 2.1  this subcomplex is the retract of  $I^{m+1}$ and 
therefore the cross section over this subcomplex induced
by  the cross section  
$s^{n}$  has  an extension  over $I^{m+1}$. Consider 
its restriction  over the $n$-skeleton  of $I^{m+1}. $  On the $m$-cube 
$(0,x_{2},x_{3}, \cdots ,x_{m+1})$ it defines an $n$-dimensional cell-cochain 
with  coefficients in  $\pi_{n}(F)=\pi.$  Obviously, this cochain is a cocycle. 
We can define  $\kappa_{z}(\sigma^{m+1}) $ to be this cocycle regarded as an $m$-cube of 
$L(\pi,n).$ 

The reasoning above suggests the definition of $\kappa_{z^{n+1}}$ 
as follows. As above, we consider 
the map $\psi : I^{m+1} @>>> \Delta_{m+1}$ and also we consider
$I^{m}$ as the face $(0,x_{2},x_{3}, \cdots ,x_{m+1})$.
The $n$-cochain $a^{m+1,-m}(\sigma^{m+1})=\lambda^{n}$ of $I^{m}$  is uniquely 
defined as follows: it is zero on  the $n$-cubes  of $I^{m}$
which have a form differing from 
$$  u=(0,\varepsilon_{1 },\varepsilon_{2},\dots,
        \varepsilon_{k},x_{i_{1}},1,1,1,
        \dots,1,x_{i_{2}},1,1,1,\dots,
        1,x_{i_{n}},\nu_{1},\nu_{2},\dots,\nu_{s}). $$
For the above $n$-cube we assume  $\varepsilon_{j}$ to be the last 
$\varepsilon$ equal to zero  in the sequence 
$\varepsilon_{1},\varepsilon_{2},\ldots,\varepsilon_{k} $  
and consider the $(n+1)$-cube in $I^{m+1}$
$$  (1, \dots,1,x_{i_{0}},1,\dots,1,x_{i_{1}},1,\dots,1,x_{i_{2}},
        1,\dots,1,x_{i_{n}},\nu_{1},\nu_{2},\dots,\nu_{s}), $$
where  $\varepsilon _{j}$ is replaced by $x_{i_{0}}$. 
If all $\varepsilon$'s are equal to 1, we consider the $(n+1)$-cube
$$  (x_{i_{0}},1,\dots,1,x_{i_{1}},
        1,\dots,1,x_{i_{2}},1,\dots,1,x_{i_{n}},
        \nu_{1},\nu_{2},\dots,\nu_{s}).         $$
This $(n+1)$-cube of $I^{m+1}$ defines via the map $\psi^{m+1}$ 
the nondegenerate 
(n+1)-simplex $v$ in the simplex $\Delta^{m+1}$ and we set  
$$  \lambda(u)=(\sigma^{m+1})^{*}(z^{n+1})(v), \;  \;  \;   
        \sigma^{m+1}: \Delta^{m+1} \rightarrow B.       $$

$\lambda^{n}$ is the cocycle. This is proved directly.

Let $\kappa_{z^{n+1}}(\sigma^{m+1})$ be the cocycle $\lambda^n$ 
regarded as the $m$-cube of $L(\pi,n)$.

It remains to show that  equalities (4.2) are fulfilled. 
Each of the four equalities (4.2) is easy to check by using 
the algebraic definition of $\kappa_{z^{n+1}}$. It is obvious that $\tilde{a}$ 
and hence $a(z^{n+1})$ is functorial. 
\end{pf}

We next define
$$      Z^{n+1}(B,\pi) \times C^{n}(B,\pi) 
        \rightarrow G^{n+1}(A^{*,-*}(B,\pi,n)).     $$

Assume that we are given a cocycle 
 $ z^{n+1}_{B \times I} \in Z^{n+1}(B \times I,\pi)$. 

By $i_{0} : B @>>>B \times I$ and $i_{1}:B @>>> B \times I $ we denote 
the imbeddings of $B$ as the lower and the upper bottom of $B \times I.$
We introduce the notation 

$$  z^{n+1}_{B}=i^{*}_{0} z^{n+1}_{B \times I} \;\;\; 
        \bar{z}^{n+1}_{B}= i^{*}_{1} z^{n+1}_{B \times I}.  $$
By virtue of Lemma 4.1
$$  a(z^{n+1}_{B})=i^{*}_{0}a(z^{n+1}_{B \times I}) \;\;\; 
        a(\bar{z}^{n+1}_{B})= i^{*}_{1}a(z^{n+1}_{B \times I}).  $$
   
In the lemma below we shall show that $a(z^{n+1}_{B})$ and 
 $ a(\bar{z}^{n+1}_{B})$  are the equivalent twisting elements of 
the  algebra $A^{*,-*}(B,\pi,n).$
 
Let us consider the standard prism construction
$$  w_{1}:C_{*}(B) @>>> C_{*+1}(B \times I)     $$
 subject to the condition
\begin{equation}
    dw_{1}-w_{1}d=i_{1}^{*}-i_{0}^{*}.
\end{equation}
To every singular simplex $\sigma^{m} \in S(B)$ this construction assigns 
the  singular  $(m+1)$-chain $w_{1}(\sigma^{m})$  which is the image of 
the main integral $(m+1)$-chain of the standard triangulation 
of $\Delta_{m} \times I$ by 
the map $\sigma^{m} \times id :\Delta_{m} \times I @>>> B \times I$.

The map
\begin{equation}
     w_{1}^{*}: A^{*,-*}(B \times I,\pi,n) @>>> A^{*,-*}(B,\pi,n) 
\end{equation}
can be defined by the composition
$$  C_{B} @>w_{1}>> C_{B \times I} @>x>> C_{*}^{\Box}(L(\pi,n)),\;\; 
        x \in A^{*,-*}(B \times I,\pi,n).       $$
From (4.3) we obtain
\begin{equation}
    w_{1}^{*}d=dw_{1}^{*}+i^{*}_{1}-i^{*}_{0}. 
\end{equation}

Let 
$$  u(z^{n+1}_{B \times I})=w_{1}^{*}(a(z^{n+1}_{B \times I})). $$
We see that $|u|=0, \;\; u \in F^{n}(A^{*,-*}(B,\pi,n)).$

 \begin{lemma}
 In the algebra $ A^{*,-*}(B,\pi,n)$ there holds the equality
$$  d_{A(B,\pi,n)}u(z^{n+1}_{B \times I})=
        a(\bar{z}^{n+1}_{B}) -a(z^{n+1}_{B})+
        u(z^{n+1}_{B \times I})a(\bar{z}^{n+1}_{B})-
        a(z^{n+1}_{B})u(z^{n+1}_{B \times I}).          $$
\end{lemma}

\begin{pf}
By Lemma 4.1 we have
$d_{A(B \times I, \pi,n)}a_{z^{n+1}_{B \times I}}=
a_{z^{n+1}_{B \times I}}a_{z^{n+1}_{B \times I}}.$

By virtue of the map $w_{1}^{*}$ and (4.5) the left side of the equality 
becomes 
$$  d( u(z^{n+1}_{B \times I}))-a(\bar{z}^{n+1}_{B})+a(z^{n+1}_{B}).  $$
The right side becomes 
\begin{equation}
    u(z^{n+1}_{B \times I})a(\bar{z}^{n+1}_{B})-
        a(z^{n+1}_{B}) u(z^{n+1}_{B \times I})
\end{equation}
which can be shown as follows. The standard triangulation of 
$\Delta_{m} \times I$ has as vertices of the lower bottom those of $\Delta_{m}$,
say, $b_{0} < b_{1} < b_{2} < \cdots b_{m}$, while the copies of the vertices 
of $\Delta_{m}$, say, $b'_{0} < b'_{1} < b'_{2} < \cdots b'_{m}$,  are 
the vertices of the upper bottom. It is assumed that $b(i) < b'(i).$  
Only  
$$  (b_{0} < b_{1} < b_{2} < \cdots b_{i}  < b'_{i} < b'_{i+1} 
        < \cdots b'_{m}),\;\;\; 
        i=0,1, \cdots m         $$  
are $(m+1)$-dimensional simplices. 

By the above definition $w_{1}(\sigma^{m})$ is the image of
$$  \sum _{0 \leq i\leq m}(-1)^{i}(b_{0} < b_{1} < b_{2} < \cdots 
        b_{i}  < b'_{i} < b'_{i+1} < \cdots < b'_{m}).      $$
Thus we see that the value of $\tilde{a} \tilde{a}$ on 
$w_{1}(\sigma^{m})$ is equal to
\begin{gather*}
    \sum _{0 \leq i\leq m}(-1)^{i}[\sum_{j \leq i}
        a'(b_{0} < b_{1} <  \cdots < b_{j})\cdot  \\
    \cdot a'(b_{j}< b_{j+1}< \cdots < b_{i}  < b'_{i} 
        < b'_{i+1} < \cdots < b'_{m})+  \\
    +\sum _{i \leq j\leq m}a'(b_{0} < b_{1} < b_{2} < \cdots 
        < b_{i}  < b'_{i} < b'_{i+1} < \cdots <  b'_{j})\cdot  \\
    \cdot a'( b'_{j} < b'_{j+1} < \cdots < b'_{m})].
\end{gather*}

Here $a'=(\sigma^{m} \times id)^{*}\tilde{a}$ and the two summands coincide 
with the two summands in (4.6). 
\end{pf}

\begin{definition}
Let $z^{n+1}_{B}\in Z^{n+1}(B,\pi)  $ and $ c^{n}_{B}\in C^{n}(B,\pi).$
On $B \times I$ consider the cocycle 
$z^{n+1}_{B \times I}=Pr^* z^{n+1}_{B}+\delta_{B \times I}c^n_{B \times I},$
where $c^n_{B \times I}$ is  the cochain $c^{n}_{B}$ imbedded in $B \times 0.$
Denote $u(c^n_B,z^{n+1}_{B})=u(z^{n+1}_{B \times I}).$
\end{definition}

 \begin{lemma}
   If $z^{n+1}_{B}\in Z^{n+1}(B,\pi)  $ and $ c^{n}_{B} \in C^{n}(B,\pi)$,
then $a(z^{n+1}_{B}+\delta_{B} c^{n}_{B})$ and $a(z^{n+1}_{B})$ are 
the equivalent twisting elements, the equivalence being given by the element 
$1+u(c^{n}_{B},z^{n+1}_{B}). $
Hence the map $\gamma: H^{n+1}(B,\pi) @>>> D(A^{*,-*}(B,\pi,n))$ defined by 
$a(z^{n+1})\in \gamma (h)$ if $ z^{n+1} \in h \in H^{n+1}(B,\pi)$
is  correct.
\end{lemma}

\noindent {\it Proof.\/}
By the assumptions of the lemma  the equality of 
Lemma 4.2  is identical to
\begin{align*}
    a(z^{n+1}_{B}+\delta_{B}c^{n}_{B}) & 
        =(1+u(c^{n}_{B},z^{n+1}_{B}))^{-1}
        a(z^{n+1}_{B})(1+u(c^{n}_{B},z^{n+1}_{B})) +  \\
    & +(1+u(c^{n}_{B},z^{n+1}_{B}))^{-1}  
        d_{A}(1+u(c^{n}_{B},z^{n+1}_{B})). \;\;\qed
\end{align*}

\section{Algebraic Model}

Recall that the bicomplexes $Y_{**}(B,\pi,n)=C_* (B,C_{*}^{\Box}(L(\pi,n)))$,
$Y_{**}(B,\pi,  \linebreak  n,G)=C_* (B,C_{*}^{\Box}(L(\pi,n),G))$ and
$Y^{**}(B,\pi,n,G)=C^* (B,C^{*}_{\Box}(L(\pi,n),G))$ are the modules 
over the bialgebra
$A^{*,-*}(B,\pi,n)=C^{*}(B,C_{*}^{\Box}(L(\pi,n)))$. 

 \begin{definition}
Let $F@>>>E@>>>B $
be a Serre fibration  with the fiber $F$ $K(\pi,n)$-space. Then 
$\pi _i (F)=0,\; i \neq n$ and $\pi_n(F)=\pi.$ Let 
$h^{n+1} \in H^{n+1}(B,\pi)$ be the obstruction class
of the fibration  and let $z^{n+1} \in h^{n+1}$.
Consider $a(z^{n+1}) \in T(A^{*,-*}(B,\pi,n))$ and the perturbed differential
\begin{equation}
\begin{gathered}
    d_{a(z^{n+1})} (y)=d_{Y}(y)+a(z^{n+1})y,   \\
    y \in Y_{**}(B,\pi,n),\;Y_{**}(B,\pi,n,G), \; Y^{**}(B,\pi,n,G). 
\end{gathered}
\end{equation}
The complex $Y_{**}(B,\pi,n)$ with this perturbed differential 
$Y_{z^{n+1}}(B,\pi,n)  \linebreak  =Y_{a(z^{n+1})}(B,\pi,n)$ is the (integral) homology 
model of the fibration. 
The complex $Y_{**}(B,\pi,n,G)$ with this perturbed differential
$Y_{z^{n+1}}(B,\pi,n,G)=Y_{a(z^{n+1})}(B,\pi,n,G)$ is the homology model  
with coefficients $G$. 
The complex $Y^{**}(B,\pi,n,G)$ with this perturbed differential
$Y_{z^{n+1}}^{**}(B,\pi,n,G)\!=Y_{a(z^{n+1})}^{**}(B,\pi,n)$ is the    
cohomology model of the fibration. 
\end{definition}

We have  
\begin{align*}
    Y_{**z^{n+1}}(B,\pi,n,G) & =Y_{**z^{n+1}}(B,\pi,n) \otimes G,  \\
    Y_{z^{n+1}}^{**}(B,\pi,n,G) & =Hom(Y_{**z^{n+1}}(B,\pi,n),G).
\end{align*}

The models of fibration depend on the choice of the cocycle in the obstruction 
class. However, they  are isomorphic complexes:
if $z^{n+1}, \; \bar{z}^{n+1}$ are two obstruction cocycles
of the same fibration,  then 
there is $c^{n} \in C^{n}(B,\pi)$
with $\delta c^{n}=\bar{z}^{n+1}-z^{n+1}$. By virtue of Lemma 4.3
$a(\bar{z}^{n+1})=(1+u(c^{n},z^{n+1}))\circ a(z^{n+1}).$ Hence the
chain map (2.9)
\begin{gather*}
    \varphi_{u(c^{n},z^{n+1})}:
        Y_{**\bar{z}^{n+1}}(B,\pi,n,G) @>>>  
        Y_{**z^{n+1}}(B,\pi,n,G), \\
    \varphi_{u(c^{n},z^{n+1})}(y)=(1+u(c^{n},z^{n+1}))(y), 
        \;\; y \in Y_{**}(B,\pi,n,G)    
\end{gather*}
is an isomorphism.

The same reasoning holds for the cohomology model.

\begin{theorem} 
Let $F@>>>E@>>>B $
be a Serre fibration with the fiber F $K(\pi,n)$-space
and let $Y_{**z^{n+1}}(B,\pi,n,G), \;
Y_{z^{n+1}}^{**}(B,\pi,n,G), z^{n+1} \in     Z^{n+1}(B,\pi_{n}(F)),$ be 
its homology (resp., cohomology) model. Then: 

\rom{(i)} There are chain and cochain maps defined uniquely up to  chain  
homotopy
\begin{gather*}
    Y_{**z^{n+1}}(B,\pi,n,G) @>>> C_{*}^{\Box}(E,G),  \\
    C^{*}_{\Box}(E,G) @>>> Y_{z^{n+1}}^{**}(B,\pi,n,G)
\end{gather*}
inducing an isomophism of homology groups.  

\rom{(ii)} if $Y_{**\bar{z}^{n+1}}(B,\pi,n,G), \;
Y_{\bar{z}^{n+1}}^{**}(B,\pi,n,G)$ are other models of the fibration
and $\delta c^{n}=\bar{z}^{n+1}-z^{n+1}, \;c^{n} \in C^{n}(B,\pi_{n}(F)),$ 
then  the triangles
\vskip+0.5cm

\unitlength=1.00mm
\special{em:linewidth 0.4pt}
\linethickness{0.4pt}
\begin{picture}(94.00,25.00)
\put(3.00,25.00){\makebox(0,0)[cc]{$C_*^{\Box}(E)$}}
\put(33.00,25.00){\makebox(0,0)[cc]{$Y_{**z^{n+1}}(B,\pi,n)$}}
\put(33.00,5.00){\makebox(0,0)[cc]{$Y_{**\bar{z}^{n+1}}(B,\pi,n),$}}
\put(33.00,10.00){\vector(0,1){10.00}}
\put(34.00,15.00){\makebox(0,0)[lc]{$\scriptstyle \varphi_{u(c^n,z^{n+1})}$}}
\put(20.00,24.67){\vector(-1,0){11.00}}
\put(20.00,10.00){\vector(-4,3){12.00}}
\put(60.00,25.00){\makebox(0,0)[cc]{$C^*_{\Box}(E,G)$}}
\put(93.00,25.00){\makebox(0,0)[cc]{$Y_{z^{n+1}}^{**}(B,\pi,n,G)$}}
\put(93.00,5.00){\makebox(0,0)[cc]{$Y_{\bar{z}^{n+1}}^{**}(B,\pi,n,G)$}}
\put(94.00,15.00){\makebox(0,0)[lc]{$\scriptstyle \varphi_{u(c^n,z^{n+1})}$}}
\put(68.33,25.00){\vector(1,0){10.67}}
\put(93.00,20.00){\vector(0,-1){10.00}}
\put(67.67,18.67){\vector(3,-2){12.33}}
\end{picture}

\noindent are commutative up to  chain homotopy.
\end{theorem}

The proof below uses the geometric interpretation of the model.

Let $K$ be a simplicial set, $z^{n+1} \in Z^{n+1}(K,\pi)$,  and define
the cubical complex $K \times_{z^{n+1}} L(\pi,n)$  as follows.
Consider as $(p+q)$-dimensional cubes of 
$K \times_{z^{n+1}} L(\pi,n)$ the pairs $(\sigma^{p},\tau^{q})$, where  
$\sigma^{p} \;$ is a $p$-dimensional simplex of $K$  
and $\tau^{q} \;$ is a $q$-dimensional cube
of the cubical set $L(\pi,n).$ 

The face operators are defined by virtue of (2.4) and (4.2) 
as follows. Let $  \kappa_{z^{n+1}}: K @>>> L(\pi,n) $ be map (4.1). 

Define
\begin{gather*}
    d^{1}_{p+i}(\sigma^{p},\tau^{q})=
        (\sigma^{p},d^{1}_{i}\tau^{q}), \;\;\; 1 \leq i \leq q,  \\
    d^{0}_{p+i}(\sigma^{p},\tau^{q})=
        (\sigma^{p},d^{0}_{i}\tau^{q}), \;\;\; 1 \leq i \leq q, \\
    d^{1}_{i}(\sigma^{p},\tau^{q})=
        (\sigma^{p}_{i-1},\tau^{q}),    \;\;\; 1 \leq i \leq p, \\
        d_{p}^{0}(\sigma^{p},\tau^{q})=(\sigma^{p}_{p},\tau^{q}), \\
    d^{0}_{i}(\sigma^{p},\tau^{q})=
        (\sigma^{p}_{1},
        \kappa_{z^{n+1}}(\sigma^{p}_{2}) \circ \tau^{q}) \;\;\;
        1\leq i\leq (p-1),
\end{gather*}
where $\sigma^{p}_{1}$ is the first $(i-1)$-face and $\sigma^{p}_{2}$ is the
last $(p-i+1)$-face of $\sigma^{p}$, while $\circ$ is the product in the ring 
$\tilde{C}_*^{\Box}(L(\pi,n))$.

The degeneracy operators are defined only partially:
$$  s^{p+i}(\sigma^{p},\tau^{q})=
        (\sigma^{p},s^{i}(\tau^{q})),\;\;\; 
        1 \leq i \leq (q+1).        $$

The chain complexes $\tilde{C}_{*}^{\Box}(K \times_{z^{n+1}} L(\pi,n))$, 
$\bar{C}_{*}^{\Box}(K \times_{z^{n+1}} L(\pi,n))$,   \linebreak 
$C_{*}^{\Box}(K \times_{z^{n+1}} L(\pi,n))$, 
$C^{*}_{\Box}(K \times_{z^{n+1}} L(\pi,n),G)$ of $K \times_{z^{n+1}} L(\pi,n))$
are defined as for the case of cubical  sets.

The obvious fact is

\begin{lemma}
The integral chain complex of the cubical 
complex   \linebreak  $S(B) \times_{z^{n+1}} L(\pi,n)$ is  the complex
$Y_{**z^{n+1}}(B,\pi,n)$.
\end{lemma}

Let $|K \times_{z^{n+1}} L(\pi,n)|$
be the Milnor realization of this complex (Section 2). 
Let $||K||$ be the Giever--Hu realization of the simplicial set $K$
(i.e., the degeneracy operators are passive).
We have the map
\begin{equation}
\alpha :    |K \times_{z^{n+1}} L(\pi,n)|   @>>> ||K||
\end{equation}
defined for every $(\sigma^{p},\tau^{q})$ by
$   I^{p+q}  @>proj.>> I^{p} @>\psi>> \Delta_{p}$.      
The complex $ |K \times_{z^{n+1}} L(\pi,n)|$ is filtered  
by  subcomplexes
$   F^{r}=\underset{p \leq r}{\cup}|(\sigma^{p},\tau^{q})|$.    
We have 

\begin{lemma}
$H_{r+q}(F^{r},F^{r-1})=C_{r}(K,H_{q}^{\Box}(L(\pi,n))).$
\end{lemma}

\begin{pf}
$F^{r}/F^{r-1}$ is a wedge of $CW$-complexes, one complex $K_{\sigma^{r}}$ 
for every $\sigma^{r} \in K$. Hence $H_{r+q}(F^{r},F^{r-1})=
\sum_{\sigma^{r}}H_{r+q}(K_{\sigma^{r}}).$ The filtration of $K_{\sigma^{r}}$
by its skeletons gives $H_{r+q}(K_{\sigma^{r}})=H_{q}^{\Box}(L(\pi,n)).$
\end{pf}

Consider the standard  map 
$   ||S(B)||  @>>> B    $ 
of the Giever--Hu realization of $S(B)$ and also consider the induced fibration
$   E' @>pr>> ||S(B)||$.    
In the diagram
$$
    \begin{CD}
        E @<<< E'  \\
        @V{pr}VV   @VV{pr}V  \\ 
        B @<<< ||S(B)||
    \end{CD}
$$
the horizontal maps induce an isomorphism of homology. In the induced fibration 
consider the filtration given by $pr^{-1}(||S(B)||^{r})$.

\begin{proposition} 
There is a commutative diagram
\vskip+0.5cm

\unitlength=1.00mm
\special{em:linewidth 0.4pt}
\linethickness{0.4pt}
\begin{picture}(65.00,25.00)
\put(25.00,25.00){\makebox(0,0)[cc]{$E'$}}
\put(45.00,25.00){\vector(-1,0){15.00}}
\put(65.00,25.00){\makebox(0,0)[cc]{$|S(B)\times_{z^{n+1}}L(\pi,n)|$}}
\put(40.34,5.00){\makebox(0,0)[cc]{$\|S(B)\|\,,$}}
\put(26.00,20.00){\vector(2,-3){6.67}}
\put(54.00,20.00){\vector(-2,-3){6.67}}
\put(53.33,15.00){\makebox(0,0)[lc]{$\scriptstyle \alpha$}}
\put(26.00,15.00){\makebox(0,0)[rc]{$\scriptstyle pr$}}
\end{picture}

\noindent where the upper map is the map of filtered spaces and induces an isomorphism 
of the first terms of the related spectral sequences.
\end{proposition}

\begin{pf}
The map $|S(B) \times_{z^{n+1}} {}_{z} L(\pi,n)| @>>> E'$  is  constructed 
by induction on degree of cells. Let $s$ be a cross section over 
the $n$-skeleton of $||S(B)||$ whose obstruction cocycle is $z^{n+1}.$
The induction steps are as follows:

(0). $(\sigma^{0},\tau^{0})$ is a vertex.  
  Let $f|(\sigma^{0},\tau^{0})=s \alpha |(\sigma^{0},\tau^{0})=s(\sigma^{0}).$ 

(i), $i \!<\! n.\;(\sigma^{j},\tau^{i-j}) $ are $i$-dimensional 
cells.\;Let $f|(\sigma^{j},\tau^{i-j})\!=\!s \alpha |(\sigma^{j},\tau^{i-j}).$   

(n). $(\sigma^{j},\tau^{n-j})$ are  $n$-cells.
Let $f|(\sigma^{j},\tau^{n-j})=s \alpha |(\sigma^{j},\tau^{n-j}), \;\; j < n$, 
and let $f|(\sigma^{0},\tau^{n})$ be the map of this ``$n$-sphere'' in a fiber
over  $\sigma^{0}$ as an element of $\pi_{n}(F,s(\sigma^{0})).$

(n+1). For the cell $(\sigma^{n+1},\tau^{0})$ the map is already defined for its 
boundary and the image lies over $\sigma^{n+1}$. This $n$-sphere is 
homotopic to $0$ over $\sigma^{n+1}$ by virtue of the fact 
that $z^{n+1}(\sigma^{n+1}) \in \pi_{n}(F)$  is the class 
of the  $n$-sphere defined by $s$ on the boundary of $\sigma^{n+1}$.
Hence $f$ extends from the boundary onto the whole cell and the image lies over
$\sigma^{n+1}.$

The map $f$  of the boundary of $(\sigma^{1},\tau^{n})$ is homotopic to $0$ and 
hence it extends onto the whole cell.

For the rest of the $(n+1)$-cells $(\sigma^{i},\tau^{n+1-i}), \;\;
1< i \leq n$, we assume that 
$f|(\sigma^{i},\tau^{n+1-i})=s \alpha |(\sigma^{i}, \tau^{n+1-i}).$

(m), $ m>n+1$.  $(\sigma^{i},\tau^{j}), \;\; i+j=m,$
$f$ already defined on its boundary is an $(m-1)$-sphere over $\sigma^{i}$ 
and by the fact that $m-1>n$ it is homotopic to $0$. $f$ extends over the whole 
cell.
The map is now constructed.  By the inductive construction we see that it
preserves filtrations.

The first term of both filtrations in the proposition under consideration 
is  $C_{p}(S(B),H_{q}(F))$, and by virtue of the above map the homomorphism 
of the first terms of the spectral sequences is an isomorphism. 
\end{pf}

\noindent {\it Proof of Theorem $5.1$.}
The standard imbedding 
$$  S(B) \times_{z^{n+1}} L(\pi,n) @>>> 
        Q(|S(B) \times_{z^{n+1}} L(\pi,n)|) $$
gives the chain map
$$  C_{*}^{\Box}(S(B) \times_{z^{n+1}} L(\pi,n)) @>>> 
        C_{*}^{\Box}(|S(B) \times_{z^{n+1}} L(\pi,n)|)  $$
inducing an isomorphism of cubical singular homologies. Then by Proposition 5.1 
the composition 
\begin{gather*}
    Y_{**z^{n+1}}(B,\pi,n,G)=C_{*}^{\Box}(S(B) 
        \times_{z^{n+1}} L(\pi,n),G) @>>>  \\
    @>>> C_{*}^{\Box}(|S(B) \times_{z^{n+1}} 
        L(\pi,n)|,G) @>f_{z^{n+1}_{*}}>> 
        C_{*}^{\Box}(E,G)           
\end{gather*}
induces an isomorphism of homology. Let this chain map be the map  of (i).

Consider the  fibration $F @>>>E'' @>>> B \times I$
induced by projection $pr:B \times I @>>> B$ and  a  map as in  part 
(i) of the theorem: 
$$  E'' @<f_{z^{n+1}_{B \times I}}<< 
        | S(B \times I) \times_{z^{n+1}_{B \times I}} L(\pi,n)|   $$
where $z^{n+1}_{B \times I}=pr*z_B^{n+1}-\delta\bar{c}^n$
$(\bar{c}^n$ is $c^n$ placed on $B\times I$) is such that
$   i_{0}^{*}z^{n+1}_{B \times I}=z^{n+1}_{B}$, 
$   i_{1}^{*}z^{n+1}_{B \times I}=\bar{z}^{n+1}_{B}$.
We  can choose this map so that  on the bottoms of $B \times  I$ it will 
coincide  with $f_{z^{n+1}}$ and $f_{\bar{z}^{n+1}}$. Consider
the imbeddings 
\begin{gather*}
    C_{*}^{\Box}( |S(B \times I) \times_{z^{n+1}_{B \times I}} L(\pi,n)|)  
        @<<< C_{*}^{\Box}( |S(B) \times_{z^{n+1}_{B}} L(\pi,n)|),  \\
    C_{*}^{\Box}( |S(B \times I) \times_{z^{n+1}_{B \times I}} L(\pi,n)|)
        @<<< C_{*}^{\Box}( |S(B) \times_{\bar{z}^{n+1}_{B}} L(\pi,n)|).
\end{gather*}
In view of the chain map  
$$      C_{*}^{\Box}(E) @<<<  C_{*}^{\Box}(E'') @<<<  
        C_{*}^{\Box}(|S(B \times I) \times_{z^{n+1}_{B \times I}} 
        L(\pi,n)|)      $$
it is enough to show that the  triangle
\vskip+0.5cm

\unitlength=1.00mm
\special{em:linewidth 0.4pt}
\linethickness{0.4pt}
\begin{picture}(82.00,24.67)
\put(23.00,24.33){\makebox(0,0)[cc]{$C_*^{\Box}(|S(B\times I)\times_{z_{B\times I}^{n+1}}L(\pi,n)|)$}}
\put(81.00,10.00){\vector(0,1){10.00}}
\put(82.00,15.00){\makebox(0,0)[lc]{$\scriptstyle \varphi_{u(c^n,z_B^{n+1})}$}}
\put(60.00,24.67){\vector(-1,0){11.00}}
\put(60.00,10.00){\vector(-4,3){12.00}}
\put(81.33,24.33){\makebox(0,0)[cc]{$C_*^{\Box}(|S(B)\times_{z_{B}^{n+1}}L(\pi,n)|)$}}
\put(81.33,5.33){\makebox(0,0)[cc]{$C_*^{\Box}(|S(B\times I)\times_{z_{B\times I}^{n+1}}L(\pi,n)|)$}}
\end{picture}

\noindent is   commutative up to  chain homotopy.
To this end consider
$$      C_{*}^{\Box}( |S(B \times I) \times_{z^{n+1}_{B \times I}} L(\pi,n)|)  
        @<\varphi_{u(\bar{c}^{n},pr^{*}z^{n+1}_{B})}<< 
    C_{*}^{\Box}( |S(B \times I) \times_{pr^{*}z^{n+1}_{B}} L(\pi,n)|).   $$
The two imbeddings of
$C_{*}^{\Box}( |S(B) \times_{z^{n+1}_{B}} L(\pi,n)|) $
on the left side complex are chain homotopic and the composition  gives 
the required homotopy. This proves (ii) for an integral homology. 
Hence we obtain (ii) for the coefficient group  $G$ and cohomology.

In particular, if $\bar{z}^{n+1}=z^{n+1}$ and $c^{n}=0$, then
we deduce that the map in (i) 
is defined uniquely up to  chain homotopy. 
\;\;\qed

\section{Multiplicative Structure in the Algebraic Model}

We have seen that the cohomology model $Y^{**}_{z^{n+1}}(B,\pi,n,G)$
of a fibration $F\to E\to B$ with $F$ a $K(\pi,n)$-space 
can be identified with the cochain complex of the cubical complex 
$ S(B) \times_{z^{n+1}} L(\pi,n) $.
Each cubical complex is endowed with multiplicative structure via 
the Serre cup product. Hence $Y^{**}_{z^{n+1}}(B,\pi,n,\Lambda),$ where
$\Lambda $ is a commutative ring, has a multiplicative structure. It 
can evidently be described as follows.

Let $x^{s}$, $\;y^{t} \in Y^{**}_{z^{n+1}}(B,\pi,n,\Lambda)$. 
Being an element of the module  $C^{p}(B,\linebreak C^{q}(L(\pi,n)\Lambda)),$
the $(p,q)$-component, where $p+q=s+t$, of the product
$w^{s+t}=x^{s}y^{t}$  is a $p$-dimensional cochain 
 of $B$ with coefficients in $C^{q}(L(\pi,n),\Lambda)$ defined
 in the following manner:  for the pair $(\sigma^{p},\tau^{q})$
consider all decompositions of $1,2, \cdots ,(p+q)$ in two disjoint 
sets $(H,K)$ and let 
$$  w^{s+t}[(\sigma^{p},\tau^{q})]=
        \sum_{(H,K)}(-1)^{a(H,K)}x(d^{0}_{K}
        (\sigma^{p},\tau^{q}))y(d^{1}_{H}(\sigma^{p},\tau^{q})).  $$

\begin{theorem}
The cochain map of Theorem $5.1$ 
$$  C^{*}_{\Box}(E,\Lambda) @>>> Y_{z^{n+1}}^{**}(B,\pi,n,\Lambda)  $$
is multiplicative.
\end{theorem}

\noindent {\it Proof.} The above cochain map is induced (see the proof of 
Theorem 5.1) by the map of cubical sets
$   S(B) \times_{z^{n+1}} L(\pi,n) @>>> Q(E). \;\; \qed     $

\section{Action of the Group $H^n(B,\pi)$ on the Homology
of the Complexes $Y_{**z^{n+1}}(B,\pi,n,G)$ and $Y^{**}_{z^{n+1}}(B,\pi,n,G)$ }

Consider the space $B \times \Delta_{2}$ and assume that we are given 
a cocycle 
$ z^{n+1}_{B \times \Delta^2} \in Z^{n+1}(B \times \Delta_{2}, \pi).$
  
We have three imbeddings
$   i_{0},i_{1},i_{2}, :    B \rightarrow B \times \Delta^{2}   $  
for three vertices of $\Delta_{2}$ and three imbeddings 
$   i_{01},i_{12},i_{02}, : B \times I \rightarrow B\times \Delta^{2}  $
for three 1-faces of $\Delta_{2}.$
For a given cocycle $ z^{n+1}_{B \times \Delta ^{2}} \in Z^{n+1}(B \times \Delta^{2})$
consider six (n+1)-cocycles (three on $B$ and three 
on $B\times I $): 
\begin{align*}
    z^{n+1}_{B,0 } =i^{*}_{0}( z^{n+1}_{B \times \Delta ^{2}} ), 
        \; \; \; \; \;\; & z^{n+1}_{B \times I ,01} 
        =i^{*}_{01}(  z^{n+1}_{B \times \Delta ^{2}} ) ,  \\
    z^{n+1}_{B ,1} =i^{*}_{1}( z^{n+1}_{B \times \Delta ^{2}} ), 
        \; \; \; \; \;\; &  z^{n+1}_{B \times I ,12} 
        =i^{*}_{12}(  z^{n+1}_{B \times \Delta ^{2}} ) ,  \\
    z^{n+1}_{B ,2} =i^{*}_{2}(  z^{n+1}_{B \times \Delta ^{2}} ), \; \; 
        \; \; \; \; &  z^{n+1}_{B \times I ,02} 
        =i^{*}_{02}(  z^{n+1}_{B \times \Delta ^{2}} ) .
\end{align*}

The left cocycles define on $B$ three twisting elements
$   a(z^{n+1}_{B ,0})$,  $a(z^{n+1}_{B ,1})$,  $a(z^{n+1}_{B ,2})$,     
while the right cocycles define on $B$ three 0-elements as 
in Lem\-ma 4.2:
\begin{gather*}
    u_{01}(z^{n+1}_{B \times I ,01}), u_{12}(z^{n+1}_{B \times I ,12}), 
        u_{02}(z^{n+1}_{B \times I ,02}) 
        \in A^{*,-*}(B,\pi,n)=   \\
    =C^{*}(B,C_{*}^{\Box}(L(\pi,n))).       
\end{gather*}

Consider the standard  map
$$      w_{2}:C_{*}(B) @>>> C_{*+2}(B \times \Delta_{2})    $$
subject to the condition
\begin{equation}
    d_{B \times \Delta_{2}}w_{2}-w_{2}d_{B}=
        (i_{01}^{*}+i_{12}^{*}-i_{02}^{*})w_{1}.
\end{equation}
To every singular simplex 
 $\sigma^{m} \in S(B)$  this map assigns the  singular  $(m+2)$-chain 
$w_{2}(\sigma^{m})$  which is the image of 
the main integral    $(m+2)$-chain of the standard triangulation 
of $\Delta_{m} \times \Delta_{2}$  by the map $\sigma^{m} \times id :
\Delta_{m} \times \Delta_{2} @>>> B \times \Delta_{2}. $ 

Define the map
$$  w^{*}_{2}: A^{*,-*}(B \times \Delta_{2},\pi,n) @>>> 
        A^{*,-*}(B,\pi,n)       $$
by the composition
$$  C_{*}(B) @>w_{2}>> C_{*}(B \times \Delta_{2}) 
        @>x>> C_{*}^{\Box}(L(\pi,n)),\;\; 
         x \in A(B \times \Delta_{2},\pi,n).    $$

Let
$   v(z^{n+1}_{B \times \Delta_{2}})=
        w^{*}_{2}(a(z^{n+1}_{B \times \Delta_{2}}))$.
We see that $|v|=-1$ and $v=v^{0,-1}+v^{1,-2}+v^{2,-3}+ \cdots $.

 \begin{lemma}
 In the above notation we have the equality
$$      d_{A(B,\pi,n)}v=u_{01 }+u_{12} +u_{01}u_{12}-u_{02}+
        v a(z^{n+1}_{B ,0})+ 
        a(z^{n+1}_{B ,2})v       $$ 
in the algebra $ A^{*,-*}(B, \pi,n).$ 
\end{lemma}

\begin{pf} 
The lemma is proved similarly to Lemma 4.2, using equality (7.1). 
\end{pf}

In what follows we abbreviate $ Y_{**z^{n+1}}(B,\pi,n) $ to  $Y_{z^{n+1}}.$

In the situation we are considering we have the triangle of 
chain isomorphisms 
\vskip+0.5cm
\unitlength=1.00mm
\special{em:linewidth 0.4pt}
\linethickness{0.4pt}
\begin{picture}(61.33,28.00)
\put(25.00,25.00){\makebox(0,0)[cc]{$Y_{z_{B,0}^{n+1}}$}}
\put(35.00,25.00){\vector(1,0){15.00}}
\put(60.00,25.00){\makebox(0,0)[cc]{$Y_{z_{B,1}^{n+1}}$}}
\put(60.00,5.00){\makebox(0,0)[cc]{$Y_{z_{B,2}^{n+1}}\;.$}}
\put(60.00,20.00){\vector(0,-1){10.00}}
\put(61.33,15.00){\makebox(0,0)[lc]{$\scriptstyle \varphi_{u_{12}}$}}
\put(32.00,20.00){\vector(2,-1){19.00}}
\put(38.33,14.00){\makebox(0,0)[rc]{$\scriptstyle \varphi_{u_{02}}$}}
\put(43.00,28.00){\makebox(0,0)[cb]{$\scriptstyle \varphi_{u_{01}}$}}
\end{picture}

 \begin{corollary}
 The above triangle is commutative up to  chain homotopy.
\end{corollary}

\begin{pf}
 By virtue of  Lemma~7.1 the 
map $F:Y_{z^{n+1}_{B,0}} \rightarrow Y_{z^{n+1}_{B,2}}$ defined 
by $F(x)=vx$
yields a chain homotopy between the composition $\varphi_{u_{12}} \varphi_{u_{01}}$
and $ \varphi_{u_{02}}. $ 
\end{pf}

 \begin{corollary}
 Let $z^{n+1}_{B} \in Z^{n+1}(B,\pi)$ and $c^{n},c^{n}_{1} \in C^{n}(B,\pi).$
Then the triangle
\vskip+0.5cm

\unitlength=1.00mm
\special{em:linewidth 0.4pt}
\linethickness{0.4pt}
\begin{picture}(76.33,28.00)
\put(29.00,25.00){\makebox(0,0)[cc]{$Y_{z^{n+1}+\delta(c^n+c_1^n)}$}}
\put(75.00,25.00){\makebox(0,0)[cc]{$Y_{z^{n+1}+\delta c^n}$}}
\put(75.00,5.00){\makebox(0,0)[cc]{$Y_{z{n+1}}$}}
\put(75.00,20.00){\vector(0,-1){10.00}}
\put(76.33,15.00){\makebox(0,0)[lc]{$\scriptstyle \varphi_{u(c^n,z^{n+1})}$}}
\put(42.00,20.00){\vector(2,-1){19.00}}
\put(48.33,14.00){\makebox(0,0)[rc]{$\scriptstyle \varphi_{u(c_1^n+c^n,z^{n+1})}$}}
\put(53.00,28.00){\makebox(0,0)[cb]{$\scriptstyle \varphi_{u(c_1^n,z^{n+1}+\delta c^n)}$}}
\put(42.00,25.00){\vector(1,0){23.00}}
\end{picture}

\noindent is commutative up to chain homotopy.  
\end{corollary}

\begin{pf}
 Let $Pr:B \times \Delta^{2} \rightarrow B $ be the projection  and
consider the cocycle of $B \times \Delta ^{2}$
$   z^{n+1}_{B \times \Delta^{2}}=Pr^{*}(z^{n+1}_{B})+\delta _{B \times
        \Delta^{2}}(\bar{c}^{n})+ 
        \delta_{B \times \Delta^{2}}(\bar{\bar{c}}^{n}+
        \bar{\bar{c}}^{n}_{1})$, 
$\bar{c}^{n}$ being $c^{n}$
identified as the cochain of the subcomplex $B \times 1$ and zero otherwise,
 and $\bar{\bar{c}}^{n}+\bar{\bar{c}}^{n}_{1}$ being  
$ c^{n}+c^{n}_{1}$   identified as the cochain of $B \times 2$ and zero 
 otherwise. The rest follows from Corollary 7.1.
\end{pf}

 \begin{corollary}
 Let $z^{n+1}_{B} \!\in \!Z^{n+1}(B,\pi)$ and $0\! \in\! C^{n}(B,\pi)$; then 
$ \varphi_{u(0,z^{n+1})}: Y_{z^{n+1}_{B}} \rightarrow Y_{z^{n+1}_{B}}$ 
is homotopic to the identity. 
\end{corollary}

\begin{pf} 
 From $0=0+0$ and Corollary 7.1 we see that 
$\varphi_{u(0,z^{n+1})} \varphi_{u(0,z^{n+1})}
\sim \varphi_{u(0,z^{n+1})};$ hence  by the fact that $\varphi_{u(0,z^{n+1})}$ 
is surjective, $\varphi_{u(0,z^{n+1})}$ is homotopic to the identity 
(in fact, $\varphi_{u(0,z^{n+1})}$ is the
identity, but we do not need  this). 
 \end{pf}

 \begin{corollary} 
 For $z^{n+1}_{B} \in Z^{n+1}(B,\pi)$ and $c^{n} \in C^{n}(B,\pi)$
the triangle
\vskip+0.5cm

\unitlength=1.00mm
\special{em:linewidth 0.4pt}
\linethickness{0.4pt}
\begin{picture}(66.33,28.00)
\put(25.00,25.00){\makebox(0,0)[cc]{$Y_{z^{n+1}}$}}
\put(65.00,25.00){\makebox(0,0)[cc]{$Y_{z^{n+1}-\delta c^n}$}}
\put(65.00,5.00){\makebox(0,0)[cc]{$Y_{z{n+1}}$}}
\put(65.00,20.00){\vector(0,-1){10.00}}
\put(66.33,15.00){\makebox(0,0)[lc]{$\scriptstyle \varphi_{u(c^n,z^{n+1})}$}}
\put(32.00,20.00){\vector(2,-1){19.00}}
\put(38.33,14.00){\makebox(0,0)[rc]{$\scriptstyle id$}}
\put(44.00,28.00){\makebox(0,0)[cb]{$\scriptstyle \varphi_{u(-c^n,z^{n+1}-\delta c^n)}$}}
\put(32.67,25.00){\vector(1,0){22.33}}
\end{picture}

\noindent is commutative up to chain homotopy.
\end{corollary}

\begin{pf}
 By virtue of Corollary~7.2  
the composition in question is homotopic to 
$ \varphi_{u(c^{n}-c^{n},z^{n+1})}= \varphi_{u(0,z^{n+1})}$. 
Corollary 7.3 accomplishes the proof.
\end{pf}

 \begin{corollary} 
 For  $z^{n+1}_{B} \in Z^{n+1}(B,\pi)$ and $c^{n-1} \in C^{n-1}(B,\pi)$
the map $\varphi_{u(\delta c^{n-1},z^{n+1})}: 
Y_{z^{n+1}_{B}} \rightarrow Y_{z^{n+1}_{B}}$ 
is homotopic to the identity.
\end{corollary}

\begin{pf}
On the complex $B \times I$ consider an $n$-cochain 
$k^{n}=\delta_{B \times I}c^{n-1}-\delta_{B}c^{n-1};$
here $c^{n-1}$ and $\delta_{B}c^{n-1}$ are 
identified as cochains in $B \times 0$ 
and zero otherwise. Then 
$\delta_{B \times I}k^{n}=\delta_{B \times I}(\delta_{B}c^{n-1}). $
Identify $k^{n}$ as a cochain of $B \times \Delta^{2} $ via 
$i_{01}:B \times I \rightarrow B \times \Delta^{2}.$ 
Consider
$   z^{n+1}_{B\times\Delta^{2}}=
        Pr^{*}z^{n+1}_{B}-\delta_{B\times\Delta^{2}}k^{n}$.   
The restrictions of 
this cocycle on $(B \times I)_{02}$ and $(B \times I)_{12}$
 are equal to $Pr^{*}(z^{n+1}_{B})$ and the restriction 
on $(B \times I)_{01}$  is
$Pr^{*}z^{n+1}-\delta_{B \times I}(\delta_{B}c^{n}).$
Hence Corollaries 7.2 and 7.4 carry the proof to the end. 
\end{pf}

 \begin{theorem}
 For every $z^{n+1} \in Z^{n+1}(B,\pi)$  and every abelian group $G$  the
cohomology group $H^{n}(B,\pi)$ acts 
from the left on the  homology and cohomology groups 
$H_{q+1}(Y_{z^{n+1}},G), \;\; H^{q+1}(Y_{z^{n+1}},G).$
The action is given by 
 the chain map 
\begin{gather*}
    \varphi_{u(c^{n},z^{n+1})}: Y_{z^{n+1}} @>>> Y_{z^{n+1}}, \;\;\;
        \varphi_{u(c^{n},z^{n+1})}(y)=(1+ u(c^{n},z^{n+1}))(y),  \\
    c^{n} \in Z^{n}(B,\pi), \; \; y \in Y_{z^{n+1}}.
\end{gather*}
\end{theorem}

\begin{pf} 
We readily prove the  theorem, using the above corollaries.
\end{pf}

\section*{Acknowledgement}

 The research described in this
    publication  was made possible in part by Grant
    No. RVC000 from the International Science Foundation.

 For the  improvement of the proof of Proposition 5.1 the
 author is  obliged to S.~Saneblidze.
\vskip+0.2cm

\centerline{\sc References}
\vskip+0.1cm

1. E. H.Brown, Twisted tensor products.
    {\em  Ann.  Math.} {\bf 69}(1959), 223--246. 

2. J. P. May,  Simplicial objects  in algebraic topology.   
    {\em Van Nostrand, Princeton,} 1967.  

3. S. MacLane,  The homology product in $K(\Pi,n)$.
    {\em Proc. Amer. Math. Soc.} {\bf 5}(1954), 642--651.

4. M. M. Postnikov,  Cubical resolutions.
    {\em Dokl. Akad. Nauk SSSR} {\bf 118}(1958), 1085--1087.

5. N. Berikashvili,  On the obstruction functor.
    {\em  Bull. Georgian Acad. Sci.} {\em$($to appear$)$.}

6. N. Berikashvili,  An algebraic model of the Postnikov construction,
    {\em Bull. Georgian Acad. Sci.} {\bf 152}(1995).

7. D. M. Kan, Abstract homotopy.
    {\em Proc. Natl. Acad. Sci. USA} {\bf 41}(1955), 1092--1096.

8. J. Milnor,  The geometric realization  of a semi-simplicial complex.
    {\em  Ann.  Math.} {\bf 65}(1957), 357--362.

9. J. B. Giever,  On the equivalence of two singular homology theories.
    {\em Ann.  Math.} {\bf 51}(1950), 178--191.

10. S. T. Hu,  On the realizability of homotopy groups and their operations.
    {\em Pacific J. Math.} {\bf 1}(1951), 583--602.

11. J. P. Serre,  Homologie singuliere des \'{e}spaces fibr\'{e}s, applications.
    {\em Ann.  Math.} {\bf 54}(1951), 429--505.

12. N. Berikashvili, On the homology theory of fibrations.
    {\em Bull. Acad. Sci. Georgia} {\bf 139}(1990),  17--19.
\vskip+0.1cm 

\centerline{(Received 21.09.1994)}

\vskip+0.1cm 

Author's address:

A. Razmadze Mathematical Institute, Georgian Academy of Sciences

1, Z. Rukhadze St., Tbilisi 380093, Republic of Georgia

\end{document}